\documentclass[12pt]{article}

\usepackage{geometry}
\usepackage{url}

\usepackage{amsmath}
\usepackage{amsfonts}
\usepackage{amsthm}
\usepackage{indentfirst}

\numberwithin{equation}{section}

\theoremstyle{plain}
\newtheorem{theorem}{\indent\rm T\,h\,e\,o\,r\,e\,m\;}[section]
\newtheorem{lemma}[theorem]{\indent\rm L\,e\,m\,m\,a\;}
\newtheorem{proposition}[theorem]{\indent\rm P\,r\,o\,p\,o\,s\,i\,t\,i\,o\,n\;}

\theoremstyle{definition}
\newtheorem{definition}[theorem]{\indent\rm D\,e\,f\,i\,n\,i\,t\,i\,o\,n\;}

\theoremstyle{remark}
\newtheorem{remark}[theorem]{\indent\rm R\,e\,m\,a\,r\,k\;}

\renewenvironment{proof}{\indent\rm P\,r\,o\,o\,f.\;}{\hfill $\square$ \\ \indent}
\makeatletter                                                  %
\renewcommand*{\@seccntformat}[1]{
  \csname the#1\endcsname\;-                                   %
}                                                              %
\renewcommand{\section}{\@startsection{section}{1}{0mm}        %
   {1.5\baselineskip}
   {1\baselineskip}
   {\indent\normalfont\normalsize\bfseries}
   }                                                           %
\renewcommand*{\@seccntformat}[1]{
  \normalfont\bfseries\csname the#1\endcsname\;-               %
}                                                              %
\renewcommand\subsection{\@startsection                        %
  {subsection}{2}{0mm}
  {1.5\baselineskip}
  {1\baselineskip}
  {\indent\normalfont\normalsize\itshape}}
\renewcommand*{\@seccntformat}[1]{
  \normalfont\bfseries\csname the#1\endcsname\;-               %
}                                                              %
\renewcommand\subsubsection{\@startsection                     %
  {subsubsection}{2}{0mm}
  {1.5\baselineskip}
  {1\baselineskip}
  {\indent\normalfont\normalsize\texttt}}
\makeatother                                                   %
\usepackage{stix}
\usepackage{enumitem}
\usepackage{enumitem}

\setlist[itemize]{leftmargin=0.7cm}
\usepackage{xspace}

\def\FF{\mathbb{F}}
\def\GG{\mathbb{G}}
\def\KK{\mathbb{K}}
\def\LL{\mathbb{L}}
\def\NN{\mathbb{N}}
\def\QQ{\mathbb{Q}}
\def\Q{\mathbb{Q}}

\def\ZZ{\mathbb{Z}}

\def\calo{\mathcal{O}}

\def\ab{\mathrm{ab}}

\def\Aut{\mathrm{Aut}}

\def\Frob{\mathrm{Frob}}
\def\Gal{\mathrm{Gal}}
\def\GL{\mathrm{GL}}

\def\PGL{\mathrm{PGL}}

\def\SL{\mathrm{SL}}

\def\tors{\mathrm{tors}}

\def\Qbar{\overline{\QQ}}
\newcommand\B{{\rm (B)}\xspace}
\newcommand\ADZ{{\rm ADZ}\xspace}

\newcommand{\mug}{{\boldsymbol\mu}}

\newtheorem{conjecture}[theorem]{\indent\rm C\,o\,n\,j\,e\,c\,t\,u\,r\,e}

\newtheorem{example}[theorem]{\indent\rm E\,x\,a\,m\,p\,l\,e}

\begin{document}
\thispagestyle{empty}

\begin{center}
{\sc\large Lea Terracini}
\end{center}
\vspace {1.1cm}

\centerline{\large{\textbf{The Bogomolov Property through Galois Representations}}}

\renewcommand{\thefootnote}{\fnsymbol{footnote}}

\footnotetext{
The author is member of the Italian Indam group GNSAGA}

\renewcommand{\thefootnote}{\arabic{footnote}}
\setcounter{footnote}{0}

\vspace{0.6cm}
\begin{center}
\begin{minipage}[t]{11cm}
\small{
\noindent \textbf{Abstract.}
 
The Bogomolov property \B for an algebraic extension of \(\QQ\) asserts the existence of a uniform positive lower bound for the absolute logarithmic Weil height outside the group of roots of unity. Originally introduced as a weakening of Northcott's property and closely related to Lehmer's conjecture, it has been established for several natural classes of infinite extensions, including abelian extensions and fields generated by torsion points of elliptic curves defined over the rationals.

Given a Galois representation of an absolute Galois group, one can associate with it the algebraic extension fixed by its kernel and ask whether this extension has property \B. This point of view allows one to reinterpret classical results, and to generalize them to other Galois representations, both of geometric and non-geometric origin.

This expository paper gives a survey of the techniques and results in this framework. We present some results on modular representations, with a particular focus on the role of local \(p\)-adic information. We explain how Sen's theorem on totally ramified \(p\)-adic Lie extensions enters the proof of new criteria for the Bogomolov property, and how these criteria apply to representations with large local image.

This contribution is based on joint work with Francesco Amoroso, Andrea Conti, and Pietro Piras.

\medskip

\noindent \textbf{Mathematics~Subject~Classification:}
11G50; 11F80.

}
\end{minipage}
\end{center}

\tableofcontents
\bigskip


\section{Introduction}
Let \(h\) denote the absolute logarithmic Weil height on the field of algebraic numbers \(\Qbar\). A celebrated conjecture due to  Lehmer predicts that, once roots of unity are excluded,  there is a universal positive lower bound for $\deg(\alpha)h(\alpha)$, as $\alpha$ varies in  \(\Qbar^\times\). Lehmer's conjecture still remains unproved in its generality; nevertheless, for some special classes of fields,  it can even be sharpened. This point of view was considered for the first time in \cite{BombieriZannier2001}, where the authors introduced the notion of fields with the \emph{Bogomolov property} (property \B for short), where the height is uniformly  bounded below outside torsion points; in other words, zero is isolated in the set of heights on the multiplicative group of the field. 

For a number field, property \(\B\) is a consequence of Northcott's theorem. For infinite extensions, however, the question becomes subtler and often reflects arithmetic features of the field, such as the local behaviour of the extension, the structure of its Galois group, or its relation with objects of geometric origin.

This note is concerned with the following perspective.  Let \(\KK\) be a number field, let \(G_\KK=\Gal(\Qbar/\KK)\), and let
\[
  \rho:G_\KK\longrightarrow \Omega
\]
be a continuous representation into a separated topological group.  The field
\[
  \KK(\rho):=\Qbar^{\ker(\rho)}
\]
is the algebraic extension of \(\KK\) cut out by the representation \(\rho\).  We say that \(\rho\) has property \B if the field \(\KK(\rho)\)  has the Bogomolov property.  This language includes classical examples, such as totally real extension, cyclotomic extensions and in general abelian extensions of number fields, as well as non-abelian extensions coming from torsion points of elliptic curves.
Indeed, one of the main motivations for our approach is an important theorem due to Habegger \cite{Habegger2013}: if \(E\) is an elliptic curve defined over \(\QQ\), then the field \(\QQ(E_{\tors})\) generated by all torsion points of \(E\) has property \(\B\). This is naturally a theorem about the Galois representation on the adelic Tate module associated with \(E\), since the field cut out by this representation is precisely the field obtained by adjoining to \(\QQ\) the coordinates of the torsion points of \(E(\overline{\QQ})\).\par
Habegger's theorem provides the starting point for the developments
surveyed here. Amoroso and Terracini
\cite{AmorosoTerracini2024} extended his approach to modular forms
under strong supersingularity and large residual image assumptions.
Piras \cite{Piras2026} subsequently extended this result to forms
with nontrivial nebentypus, taking into account the role of inner twists.
A different development, due to Conti and Terracini
\cite{ContiTerracini2025}, uses Sen's theorem on ramification in
$p$-adic Lie extensions to obtain criteria for the Bogomolov property
for $p$-adic representations with large local image, without requiring
a geometric origin. More recently, Conti, Piras and Terracini
\cite{ContiPirasTerracini2026} developed an approach for modular forms
that are supercuspidal at a prime $p$, providing another setting in
which the prime-to-$p$ part of the adelic representation can be treated.

The paper is organized as follows. Section~\ref{sect:Bogomolov} recalls the basic definitions and examples concerning the Bogomolov property, together with the ADZ criterion and its interpretation in terms of Galois representations. Section~\ref{sect:modular} is devoted to modular Galois representations and to the known results toward property \B in this setting. In Section~\ref{sect:proof} we describe the strategy of the proof in the strong supersingular case and discuss two directions in which the argument can be extended, separately treating the \(p\)-part and the prime-to-\(p\) part of the adelic representations. 
Section~\ref{sect:applications} concerns \(p\)-adic representations with large local image and illustrates the scope of the main criterion through a selection of applications.

This note is expository.  It records the main definitions and statements presented in the talk, and it sketches some main tools used in the proof mechanisms.  Full proofs and more details can be found in the cited papers, especially \cite{Habegger2013,AmorosoTerracini2024,ContiTerracini2025}.

\section{The Bogomolov property of algebraic fields}\label{sect:Bogomolov}
\subsection{Weil height}
Let \(\alpha\in\Qbar\), and let \(\KK\) be a number field containing \(\alpha\).  The absolute logarithmic Weil height is defined as
\[
  h(\alpha)=\frac{1}{[\KK:\QQ]}\sum_{v}[\KK_v:\QQ_v] \log^+ |\alpha|_v,
  \qquad \log^+(x)=\log\max\{x,1\},
\]
where \(v\) runs over all places of \(\KK\), the absolute values are normalized to satisfy the product formula, and $\KK_v$ denotes the completion of $\KK$ at a place $v$.  This definition is independent of the chosen field \(\KK\) containing $\alpha$.

The following  elementary properties are easily proved:
\begin{itemize}
  \item \(h(\alpha\beta)\leq h(\alpha)+h(\beta)\) and $h(\zeta\alpha)=h(\alpha)$ if $\zeta $ is a root of unity.
  \item \(h(\alpha+\beta)\leq h(\alpha)+h(\beta)+\log 2\),
  \item \(h(\alpha^n)=|n|h(\alpha)\qquad(n\in\QQ)\),
  \footnote{For $\alpha\neq 0$ and $n=a/b\in\QQ$, with $a,b\in\ZZ$ and
$b\geq 1$, the notation $\alpha^n$ denotes any solution $\beta$
of $\beta^b=\alpha^a$. Any two such choices differ by a root
of unity and hence have the same height. }
  \item \(h(\sigma(\alpha))=h(\alpha)\qquad(\sigma\in G_\Q)\).
\end{itemize}

Moreover, the following results are well known:
\begin{itemize}
    \item \emph{Kronecker's theorem}: \(h(\alpha)=0\) for \(\alpha\in\Qbar^\times\) if and only if \(\alpha\) is a root of unity.
    \item \emph{Northcott's theorem}: there are only finitely many algebraic numbers of bounded degree and bounded height.
\end{itemize}

A classical conjecture, proposed by Lehmer in 1933, predicts that the product of the degree of a non-torsion algebraic number and its height is uniformly bounded from below.

\begin{conjecture}[Lehmer]
There exists a universal constant \(c>0\) such that
\[
  [\QQ(\alpha):\QQ]\,h(\alpha)\geq c
\]
for every \(\alpha\in \overline{\QQ}^{\times}\setminus\mug_\infty\).
\end{conjecture}

Although Lehmer's conjecture remains open in its full generality, certain fields of algebraic numbers are known to satisfy a stronger property: the height is uniformly bounded from below outside torsion.

\begin{definition}
Let \(\LL\subseteq\Qbar\) be a field.  We say that \(\LL\) has the \emph{Bogomolov property}, or property \B, if there exists a constant \(c_\LL>0\) such that
\[
  h(\alpha)\geq c_\LL
\]
for every \(\alpha\in \LL^\times\setminus\mug_\infty\).
\end{definition}

 \subsection{Fields with property \B. }\label{subsect:campiconB}

It has often been observed that most of the known examples of fields with property \B fall into two broad families.
\paragraph{Fields satisfying some local boundedness conditions.}
By Northcott's theorem, number fields have property \B. A non trivial example is the compositum $\Q^{\rm tr}$ of all totally real fields (the original proof was obtained by Schinzel, applying \cite[Corollary $1^\prime$]{Schinzel1974} to the linear polynomial $P(z)=z-\alpha$ to get the optimal lower bound). In~\cite{BombieriZannier2001} Bombieri and Zannier also show that totally $p$-adic fields, or more generally fields with bounded local degrees at some finite place have Property \B. By a result of Checcoli~\cite{Checcoli2013}, a Galois extension has uniformly bounded local degrees at {\sl every} finite place if and only if its Galois group has finite exponent. Thus (infinite) Galois extensions with Galois group of bounded exponent have Property \B. 

\paragraph{Fields obtained by adding torsion to some base field.} 
   Zannier conjectured that the maximal abelian extension $\QQ^{\ab}$ of the rational field also satisfies property \B, and this was proved in \cite{AmorosoDvornicich2000}; by Kronecker-Weber theorem $\Q^{\ab}$ is the field obtained by adjoining all torsion points of the multiplicative group. By a special case of \cite[Theorem 1.1]{AmorosoZannier2000}, property \B also holds for $\KK^{\ab}$, for any number field $\KK$ (notice however, that $\KK^\ab$ is not in general obtained by adding torsion points of some algebraic group). 
   An important class of non abelian extensions with Property \B has been exhibited by Habegger in \cite{Habegger2013}: let $E$ be an elliptic curve defined over \(\QQ\); then the field $\Q(E_\tors)$ generated by the coordinates of all torsion points of $E$ has the Bogomolov property.

\subsection{Property \B under the compositum of fields and \ADZ fields.}\label{subsect:compositum}
Property \B is not preserved under taking composita, even in the basic case of the compositum of an infinite extension and a finite one, as the following example, taken from \cite[Example 5.3]{AmorosoDavidZannier2014}, shows.

\begin{example}
Let \(\QQ^{\mathrm{tr}}\) be the maximal totally real extension of \(\QQ\). As already mentioned in Section \ref{subsect:campiconB}, it has property \B\ by \cite{Schinzel1974}. However, its compositum with the quadratic field \(\QQ(i)\) is the maximal totally CM extension \(\QQ^{\mathrm{tr}}(i)\) of \(\QQ\), which contains a sequence of elements whose heights tend to zero:
\[
\alpha_k=\left(\frac{2+i}{2-i}\right)^{1/k}.
\]
\end{example}

We recall a result of Amoroso, David and Zannier \cite[Lemma 1.5]{AmorosoDavidZannier2014}. This is a very useful criterion, since it can be used to establish property B for infinite fields which are neither abelian nor of bounded local degree at a fixed place, but rather, in some sense, cut across the two families described above.

\begin{definition}
Let $\KK$ be a number field, $v$ a non-archimedean place of $\KK$, and 
$\LL/\KK$ a Galois extension
with Galois group $\GG$. We say that the extension $\LL/\KK$ is ADZ at $v$ if
$\LL^{Z(\GG)}/\KK$
has  local degrees uniformly bounded over $v$. 
\end{definition}  
\noindent ($Z(\GG)$ denotes the center of the group $\GG). $
 
 \begin{lemma}[\ADZ]\label{lem:ADZ}
\cite[Lemma 1.5]{AmorosoDavidZannier2014} \ADZ fields have property \B.
\end{lemma}

An interesting feature of ADZ fields is that they are stable under composition with finite extensions. More generally, let \((\LL_i/\KK)\) be a family of extensions that are all ADZ at the same place \(v\). Assume that the local degrees at \(v\) of the fixed fields by \(Z(\Gal(L_i/K))\) are bounded independently of \(i\). Then the compositum \(\bigcup_i L_i\) is again ADZ at $v$ \cite[Lemma 5.7]{Piras2026}. In particular, the compositum of two extensions of $\KK$ that are both ADZ at $v$ is itself ADZ at $v$.

\subsection{Property \B\  for Galois representations} 

Habegger's theorem initiated the study of the Bogomolov property in non-abelian Galois extensions arising from arithmetic geometry; it will be  the starting point for our exposition.
 
\begin{theorem}[Habegger 2013, \cite{Habegger2013}] Let $E$ be a non CM elliptic curve defined over $\QQ$.  Then the field $\QQ(E_{\mathrm{tors}})$ has property \B.
\end{theorem}

If $E$ has complex multiplication then $\Q(E_\tors)$ is an abelian extension of a quadratic number field and thus it has the property \B by the quoted result of Amoroso and Zannier on $\KK^\ab$ \cite{AmorosoZannier2000}. However for a non-CM elliptic curve there are no number fields $\KK$ such that $\QQ(E_\tors)$ is contained in the maximal abelian extension of $\KK$ (cf. the discussion in the Introduction of \cite{Habegger2013}). 
\\
Notice that $$\QQ(E_{\mathrm{tors}})=\overline{\QQ}^{\ker{\widehat\rho_E}},$$
where
$$\widehat\rho_E:G_\QQ\longrightarrow  \mathrm{Aut}(\widehat T(E))\simeq \GL_2(\widehat \ZZ)$$
is the representation given by the action of $G_\QQ$ over the adelic Tate module associated to $E$.
The representation $\widehat\rho_E$  
is the product of all the $\ell$-adic representations
$$
\rho_{E,\ell}: \Gal(\Qbar/\Q) \longrightarrow Aut(T_\ell(E))=\GL_2(\ZZ_\ell)
$$
where  $T_\ell(E)$ is the $\ell$-adic Tate module. Indeed, by Galois correspondence $\Q(E[\ell^n])$ is the subfield of $\Qbar$ fixed by the kernel of the action of $G_\Q$ on the $\ell^n$-torsion points of $E$; moving to the projective limit we see that $\Q(E[\ell^\infty])=\QQ(\rho_{E,\ell})$. Since $\Q(E_{\tors})$ is the compositum of $\Q(E[\ell^\infty])$,   for every prime $\ell$, then it coincides with the  field fixed  by $\bigcap_\ell\ker(\rho_{E,\ell})=\ker(\widehat\rho_{E})$.\\ 
This point of view naturally leads to the following generalization. For any field $\KK\subseteq \Qbar$, we introduce the following 
\begin{definition}
Let $\Omega$ be a separated topological group, $\KK$ be a subfield of $\Qbar$ and $G_\KK=\Gal(\Qbar/\KK)$. We say that a continuous homomorphism \[\rho:G_\KK\longrightarrow \Omega\] has \emph{property \B} if the field cut out by \(\rho\),
\[
\KK(\rho):=\Qbar^{\ker(\rho)},
\]
has property \(\mathrm{(B)}\).
\end{definition}

We emphasize that this definition offers more than a change of
terminology. Since $G_\KK$ is compact and $\Omega$ is Hausdorff,
the representation $\rho$ allows us to identify the Galois group
$\Gal(\KK(\rho)/\KK)$ with a closed subgroup of $\Omega$.
The structure of this topological group may be reflected in the
arithmetic of the field $\KK(\rho)$ (for example in the decomposition
of primes), allowing a deeper understanding of the behaviour
of the Weil height on this field.

In this perspective,  the main result  of \cite{AmorosoDvornicich2000} implies (by Galois theory) that  every homomorphism of $G_\Q$ taking values in an abelian group has Property \B; and Habegger's theorem establishes property \B for the representation $\widehat\rho_E$ when $E$ is an elliptic curve defined over $\QQ$. 

We shall be particularly interested in the case where \(\rho\) is a Galois representation, that is, when \(\Omega=\GL_n(A)\), where \(A\) is a separated topological ring.

As an immediate consequence of the ADZ lemma, particularly simple to state in the language of Galois representations, we obtain the following
 
\begin{lemma}[ADZ, representation variant]\label{lem:ADZrepr}
Let \(\KK\) be a number field, and let
\[
\rho:G_\KK \longrightarrow \GL_d(A)
\]
be a Galois representation. If the image of a decomposition group \(G_v\) in the projective linear group \(\PGL_d(A)\) is finite for some non-Archimedean place \(v\) of \(\KK\), then \(\rho\) has property \(\mathrm{(B)}\).
\end{lemma}

\begin{proof}
As explained above, we may identify
$\GG=\Gal(\KK(\rho)/\KK)$ with the closed subgroup
$\rho(G_\KK)$ of $\GL_d(A)$.
Let $S$ be the subgroup of scalar matrices in $\GG$.
Then $S\subseteq Z(\GG)$, and hence
\(
\KK(\rho)^{Z(\GG)}\subseteq \KK(\rho)^S
\).
The latter is the field cut out by the projective representation
associated with $\rho$. By assumption, its decomposition group
corresponding to $G_v$ is finite. Since $\KK(\rho)^S/\KK$
is Galois, all its decomposition groups above $v$ are conjugate
and hence have the same finite order.
Thus $\KK(\rho)^S/\KK$ has uniformly bounded local degrees
at $v$, and so does its subextension
$\KK(\rho)^{Z(\GG)}/\KK$.
Therefore $\KK(\rho)/\KK$ is ADZ at $v$, and the conclusion
follows from Lemma~\ref{lem:ADZ}.
\end{proof}

Notice that, under the assumptions of Lemma~\ref{lem:ADZrepr},
the image $\rho(G_v)$ of the decomposition group may still be
infinite. Thus finiteness of the local projective image does
not imply that $\KK(\rho)/\KK$ itself has bounded local
degrees at $v$.

\medskip

 One of the main purposes of this article is to re-interpret some known results on property \B{} in terms of Galois  representations, with the aim of providing a unifying framework for a part of the theory.
\section{Modular Galois representations}\label{sect:modular} 
\subsection{Modular forms and Galois representations} Let $f=\sum_{n}a_nq^n\in S_k(\Gamma_0(N),\chi)$  be a normalized eigenform for the Hecke algebra, of weight $k\geq 2$, level $N$ and nebentypus $\chi$; 
let $\KK_f$ be the Hecke field of $f$, that is, the number field   generated by the coefficients $a_n$; and  $\calo_f$ its integer ring.
 A construction due to Shimura and Deligne attaches to $f$ a family of compatible Galois representations 
$$
\rho_{f,v}:G_\Q \longrightarrow \GL_2(\KK_{f,v})
$$
indexed by the finite places $v$ of $\KK_f$.\par
Each $\rho_{f,v}$ is characterized by the following properties
\begin{itemize}
    \item $\rho_{f,v}$ is continuous and unramified at primes $p$ not dividing $N\ell$, where $\ell$ is the rational prime which $v$ divides;
    \item for each prime $p$ not dividing $N\ell$, the characteristic polynomial of $\rho_{f,v}(Frob_p)$ is 
    \begin{equation}\label{eq:polcar} X^2-a_pX+\chi(p)p^{k-1}.\end{equation}
\end{itemize}

Since $G_\QQ$ is compact, each representation admits
a $G_\QQ$-stable $\calo_{f,v}$-lattice. Choosing such a lattice
and a basis, we may therefore regard $\rho_{f,v}$ as taking
values in $\GL_2(\calo_{f,v})$.
We fix these choices throughout; they do not affect the kernel
of the representation or the field it cuts out, but the residual
representation and the congruence filtration may depend on the
chosen lattice. The semisimplification of the residual representation,
however, is independent of this choice up to isomorphism.

For a rational prime $p$, we let 
$$
\rho_{f,p}:G_\Q \longrightarrow \prod_{v|p} \GL_2(\calo_{f,v} ).
$$
Moreover, let 
$$
\widehat\rho_f=\prod_{v}\rho_{f,v}:G_\Q\longrightarrow \GL_2(\widehat\calo_f),
$$
where $\widehat\calo_f$ is the profinite completion of $\calo_f$.

Motivated by Habegger's theorem and the results obtained so far, the following conjecture is ventured in \cite{AmorosoTerracini2024} 

\begin{conjecture}\label{conj:modularB}
Let $f\in S_k(\Gamma_1(N))$ be a normalized eigenform. Then,  
\begin{itemize}
    \item[a)]  for every non-Archimedean place $v$, $\rho_{f,v}$ has property \B;
    \item[b)] $\widehat\rho_{f}$ has property \B.
    \end{itemize}
\end{conjecture}

Notice that \(\QQ(\widehat\rho_f)\) is the compositum of the fields
\[
\QQ(\widehat\rho_f)=\bigcup_v \QQ(\rho_{f,v}),
\]
so that part~b) of Conjecture~\ref{conj:modularB} implies part~a). However, since property \B is not preserved under composita of fields, the converse is not immediate.

\paragraph{Abelian varieties of type $\GL_2$.}
The following conjecture, often attributed to S.~David, may be viewed as a natural extension of Habegger's theorem: 
\begin{conjecture}[S. David]\label{conj:David} Let $A$ be an abelian variety defined over a number field $\KK$. Then $\KK(A_{\rm tors})$ has property \B.\end{conjecture}

Conjectures \ref{conj:modularB} and \ref{conj:David} overlap to some extent. 
Indeed, when \(f\) is a cuspidal modular eigenform of weight \(2\), the representation \(\widehat{\rho}_f\) is the Galois representation of \(G_{\QQ}\) on the adelic Tate module \(\widehat{\mathrm T}(A_f)\), where \(A_f\) is an abelian variety defined over \(\QQ\) of dimension \([\KK_f:\QQ]\), equipped with an action of \(\calo_f\) \cite{Shimura1973}. 
Such abelian varieties are called \emph{of type \(\GL_2\)}, after \cite{Ribet1997}.  

Conversely, as a consequence of Serre's modularity conjecture, proved by Khare and Wintenberger in 2009 \cite[Corollary 10.2]{KhareWintenberger2009}, every abelian variety of \(\GL_2\)-type over \(\QQ\) is modular.

Therefore, Conjecture \ref{conj:modularB} predicts, in particular, property \(\mathrm{(B)}\) for abelian varieties of \(\GL_2\)-type defined over \(\QQ\).

\subsection{Property \B\ for modular representations}
Despite the generality of Conjecture \ref{conj:modularB}, the available results proving property \B for adelic representations attached to modular forms hold only under strong restrictions.

By Habegger's theorem, \B holds for modular forms of weight 2 and Hecke field equal to $\QQ$, because the abelian varieties of type $\GL_2$ they give rise are elliptic curves in this case.\\
Some examples in the case where the Hecke field has degree greater than one over $\QQ$ are provided by the following 
  \begin{theorem}\cite[Theorem 1.4] {AmorosoTerracini2024}\label{teo:AmorosoTerracini}
 Let $f=\sum_n a_nq^n \in S_k(\Gamma_0(N))$ be a normalized eigenform.\\ Assume that there exists a prime $p$ satisfying the following conditions
 \begin{description}
     \item[(P0)] $p\not | N$ and $p\geq 2k-1$;
     \item[(P1)] $a_p=0$;
     \item[(P2)] the reduction $\overline{\rho}_{f,p}: G_\QQ\to \GL_2(\calo_f\otimes \FF_p)$ has image 
    $$\widetilde G=\{\alpha\in \GL_2(\calo_f\otimes \FF_p)\ |\ \det\alpha\in(\FF_p^\times)^{k-1}\}.$$ 
 \end{description}
     Then
   $\widehat\rho_f$  has property \B.
  \end{theorem}

\paragraph{About hypotheses \textbf{(P0)}, \textbf{(P1)}, \textbf{(P2)}.}
If \(f\in S_2(\Gamma_0(N))\) has rational Hecke eigenvalues, then it falls within the scope of Habegger's theorem. The latter is unconditional, since for a non-CM elliptic curve defined over \(\QQ\) the existence of a prime \(p\) satisfying hypotheses \textbf{(P0)},  \textbf{(P1)},  and \textbf{(P2)} follows from the combination of two fundamental results.
\begin{itemize} 
    \item \emph{Serre's big image theorem} \cite[Théorème 2]{Serre1972}: Let \(\KK\) be a number field and let \(E_{/\KK}\) be an elliptic curve without complex multiplication. Then the mod-\(p\) Galois representation
\[
\rho_{E,p}:G_\KK \longrightarrow
\Aut(E[p])\simeq \GL_2(\FF_p)
\]
is surjective for all but finitely many primes \(p\).
    \item \emph{Elkies' theorem on supersingular primes \cite[Theorem 1]{Elkies1987} }: For every elliptic curve $E_{/\QQ}$, there exist infinitely many supersingular primes.
\end{itemize}

    For a non-CM modular form $f\in S_k(\Gamma_0(N))$ with no
nontrivial inner twists, an analogue of Serre's theorem,
due to Ribet \cite[Theorem 3.1]{Ribet1985}, ensures that
condition \textbf{(P2)} holds for every sufficiently large
prime $p$.  In the case of trivial nebentypus (the one we are dealing with)
it is easy to see that inner twists must arise from a quadratic Dirichlet character;
and Ribet proved that they do not appear when $N$ is square free \cite[3.9]{Ribet1980}. Therefore it is easy to produce examples of modular eigenforms with trivial nebentypus for which  property (P2) holds for all but finitely  many $p$.\\ However, the analogue of Elkies' theorem, namely the existence of infinitely many primes satisfying the strong supersingular condition $a_p =0$, is far from being established. 

On the contrary, for non-CM modular forms of weight strictly greater than \(2\), this condition appears to be rather exceptional. Evidence for this comes from two complementary directions. On the one hand, Calegari and Sardari \cite{CalegariSardari} proved that, for a fixed prime \(p\), there are only finitely many eigenforms \(f\) of a given level prime to $p$ satisfying \(a_p(f)=0\). On the other hand, for a fixed eigenform \(f\), the Atkin--Serre conjecture \cite[(4.11$_k$)]{Serre1967} predicts that \(a_p(f)=0\) for only finitely many primes \(p\). Thus, whether one fixes \(p,N\) and varies \(f\), or fixes \(f\) and varies \(p\), the vanishing of \(a_p(f)\) is expected to occur only rarely.\\
This phenomenon  is the main factor limiting the range of applicability of 
 Theorem \ref{teo:AmorosoTerracini}.

 \paragraph{Further developments.}
 A generalization of Theorem \ref{teo:AmorosoTerracini} to modular forms with nontrivial nebentypus
has been obtained by Piras in \cite{Piras2026}.
In this setting, inner twists are unavoidable: for a non-CM
form with nebentypus $\chi\neq 1$, complex conjugation on the
coefficient field gives a nontrivial inner twist associated
with $\chi^{-1}$. Although inner twists may also occur for
forms with trivial nebentypus, their necessary presence here
prevents condition \textbf{(P2)} from holding in its present form.
Thus both the statement and the proof must be reformulated
by restricting the representation to the absolute Galois group
of a suitable number field determined by the inner twists.

A recent result establishes property \B\ for a class of modular forms that are supercuspidal at a prime \(p\) \cite[Theorem A]{ContiPirasTerracini2026}. A brief account of this theorem will be given in Section~\ref{sect:supercuspidal}.

\section{The proof strategy and further progress}\label{sect:proof}
\subsection{Idea of the proof}
Although property \B\ is not preserved under taking composita, it is useful, for the purpose of analyzing the proof of Theorem \ref{teo:AmorosoTerracini} ---which closely follows Habegger's argument---to distinguish two separate components once a prime \(p\) satisfying the hypotheses of the theorem has been fixed:
\begin{itemize}
\item the \emph{prime-to-$p$ part,} concerning $\widehat\rho_f^{(p)}:G_\QQ\rightarrow \GL_2(\prod_{\ell\not=p}\calo_f\otimes \ZZ _\ell)$.
\item the \emph{ $p$-part,} concerning $\rho_{f,p}:G_\QQ\rightarrow \GL_2(\calo_f\otimes\ZZ_p)$.
\end{itemize}

Inspired by Habegger's argument,  a criterion ensuring the Bogomolov property for the product of these two representations has been developed \cite[Proposition 3.4]{AmorosoTerracini2024}. A further generalization was subsequently obtained in \cite[Theorem 4.2]{ContiPirasTerracini2026}.

\paragraph{The prime-to-$p$ part.}
The ADZ Lemma \ref{lem:ADZrepr} is enough to prove Property \B for the prime-to-$p$ part. 
Indeed, the proof can be outlined as follows: 
$$\begin{aligned} p\not | N \hbox{ and } a_p=0  &  \Longleftrightarrow \widehat\rho^{(p)} \hbox{ unramified and } \mathrm{Trace}(\mathrm{Frob}_p)=0 \\
& \buildrel{\eqref{eq:polcar}}\over\Longrightarrow \widehat\rho^{(p)} \hbox{ unramified and } \mathrm{Frob}_p^2 \hbox{ is central} \\
& \buildrel{\mathrm{(ADZ)}} \over\Longrightarrow \hbox{property } \B \hbox{ for the prime-to-$p$ part.}
   \end{aligned} $$

\paragraph{The $p$-part.}
This is the most delicate part of the proof; it is based on the study of the $p$-adic representation $\rho_p|_{G_p}$, where $G_p$ is a decomposition group at $p$. We only provide a brief overview of the key steps of the proof; the interested reader is referred to \cite{Habegger2013, AmorosoTerracini2024, Piras2026} for the details.
\begin{itemize}
\item Since $p$ does not divide the level $N$, the representation
$\rho_p|_{G_p}$ is crystalline. By
\cite[Theorem 1.2.4(ii)]{Scholl1990},
the characteristic polynomial of the crystalline Frobenius
on the associated filtered $\varphi$-module, defined using
the contravariant functor $D_{\mathrm{cris}}$, is the Hecke polynomial
\eqref{eq:polcar}, in accordance with the local Langlands
correspondence.
\item  A description due to Breuil \cite[Proposition 3.2]{Breuil2003} shows then that $\rho_p|_{G_p}$ is, up to a twist by an unramified character, induced by a power of a Lubin-Tate character: $G_{\QQ_{p^2}}\to \ZZ_{p^2}^\times$.
\item One may then appeal to Lubin-Tate theory (in particular \cite[Proposition 6.1]{Neukirch1999}) to control the higher ramification groups of the local extension $L=\QQ_p(\rho_p|_{G_p})$.
\item This yields, for every
$\alpha\in L^\times\setminus\QQ_{p^2}^\times$,
a \emph{local metric inequality:}
\[
0<|\sigma(\alpha)-\alpha|_p
<C\max\{|\alpha|_p^2,1\}, \qquad 0<C<1,
\]
for a suitable $\sigma\in G_p$, depending on $\alpha$,
and a constant $C$ independent of $\alpha$.

\item The next step is to pass from the local setting to the global one by propagating the local metric inequality to a positive proportion of the primes \(v\mid p\). This is the point at which the big image hypothesis \textbf{(P2)} on \(\bar{\rho}\) becomes essential.
    \item Then the product formula yields a lower bound for the height outside a small exceptional set. Dealing with the latter requires an additional descent argument, which we shall not discuss here.
    \end{itemize}
    
    \subsection{Working on the $p$-part: $p$-adic Lie groups and Sen's Theorem}
As discussed above, the condition \textbf(P1) considerably limits the scope of Theorem \ref{teo:AmorosoTerracini}. In order to verify Conjecture \ref{conj:modularB}, one is naturally led to seek weaker substitutes for the condition \(a_p=0\), replacing it with assumptions that are expected to hold more frequently.

If we examine the outline of the proof given in the previous section, we see that condition \textbf{(P1)} enters the argument in two distinct contexts:

\begin{itemize}
\item In the \emph{prime-to-$p$ part}  in order to know that a power of Frobenius is central and apply the ADZ lemma;
\item In the \emph{$p$-part} in order to apply Lubin-Tate theory, control ramification and produce the local metric inequality.
\end{itemize}

Although the restriction to $G_p$ of a \(p\)-adic Galois representation is usually the most challenging aspect to analyze, it is there that the most substantial progress has been made in relaxing the strong supersingularity condition.

\paragraph{A problem of Serre.} In the final part of his paper \cite{Serre1967}, Serre posed the following problem.
Let $K/\QQ_p$ be a finite extension, 
$\calo$ be a $p$-adic ring, and let
\[
\rho:G_K\longrightarrow \GL_d(\calo)
\]
be a continuous homomorphism.
Put $L=K(\rho)$,  $G=\Gal(L/K)\simeq\mathrm{im}(\rho)$.\\
\medskip
Assume that $L/K$ is totally ramified.

There are two natural filtrations on $G$:

\begin{itemize}
\item the \emph{Lie filtration}, coming from the codomain of $\rho$,
\[ \ldots G_{n+1}   \subseteq G_n\subseteq G_{n-1}\subseteq \ldots\subseteq G_0=G\]
where \[G_n=G\cap(1+p^n\mathrm{M}_d(\calo)).\]
\item  the \emph{higher ramification groups filtration}  (upper notation), coming from the domain of $\rho$,
\[ \ldots G^{n+1}   \subseteq G^n\subseteq G^{n-1}\subseteq \ldots\subseteq G^0=G\]
\end{itemize}

The Lie filtration measures how close an element of $G$
is to the identity in the $p$-adic topology of the matrix group.
The ramification filtration measures how close an automorphism
is to fixing every integral element of the extension, in terms
of the $p$-adic valuation of the difference between an element
and its image. Comparing the two filtrations therefore
translates congruence conditions on matrices into local
metric inequalities, and conversely.

Serre asked how these two filtrations are related and, in particular, how the corresponding valuations on \(G\) depend on each other and on the ramification index $e_K$ of the base field.

The same question can be formulated for any totally ramified \emph{\(p\)-adic Lie extension} \(L/K\), that is, an extension whose Galois group \(\Gal(L/K)\) is a \(p\)-adic Lie group. Indeed, the \(p\)-adic Lie structure endows \(G:=\Gal(L/K)\) with a \emph{Lie filtration} $(G_n)$ in the sense of \cite[\S 3]{Sen1973}, which can then be compared with the ramification filtration arising from the arithmetic structure of the extension.

  \paragraph{Sen's theorem}
  The answer to Serre's question was provided in \cite[\S 4]{Sen1972}\begin{theorem}[Sen]\label{teo:sen}
Let $F$ be a local field of characteristic $0$, residue
characteristic $p$, and absolute ramification index $e$.
Let $L/F$ be a totally ramified Galois extension.
Assume that $G=\Gal(L/F)$ is a $p$-adic Lie group
with $\dim(G)>0$, and let $(G_n)$ be a Lie filtration of $G$.
Then there exists a constant $c>0$ such that
\[
G^{ne+c}\subseteq G_n\subseteq G^{ne-c}
\]
for all $n\geq 0$, with $G^r=G$ for $r<0$.
\end{theorem}
\begin{remark} Although Sen's theorem was originally formulated for extensions of local fields, the requirement that the residue field is finite does not enter in a substantial way in the proof. In fact, as shown in \cite{Probst2008}, the statement remains valid for totally ramified Galois extensions of fields of characteristic \(0\), complete with respect to a discrete non-Archimedean valuation and with perfect residue field of characteristic \(p>0\). 
In particular, it holds when $F$ is the completion of an algebraic extension of $\QQ_p$ with a finite ramification index.
\end{remark}

Sen's theorem says that in a $p$-adic Lie extension the layers of a Lie filtration control the upper indices of the higher ramification groups. By Herbrand theorem, one can then pass to the lower indices and provide local metric inequalities as for the Lubin-Tate extensions (which in fact are a particular case). Details can be found in \cite[\S1]{ContiTerracini2025} for the case of finite residue field, and in \cite[\S2]{ContiPirasTerracini2026} in the general case.

It turns out that Sen's theorem provides estimates on the
ramification of totally ramified $p$-adic Lie extensions,
without requiring the explicit description furnished by
Lubin-Tate theory, which was used in the proof of the
$p$-part of Theorem~\ref{teo:AmorosoTerracini} under the
assumption $a_p=0$.

\paragraph{Property \B\ for $p$-adic Lie extensions.}
Using Sen's theorem on the local part, as outlined above, has led to criteria ensuring the Bogomolov property for \(p\)-adic Lie extensions of number fields. 

We present here the main result of \cite{ContiTerracini2025}.
Let $\KK$ be a number field,
$\LL/\KK$ be a Galois extension with $\GG=\Gal(\LL/\KK)$ a (compact) $p$-adic Lie group.
  Fix a non-Archimedean place $v$ of $\KK$ dividing $p$;  let $G=G_v$ be a decomposition group of $\GG$ at  $v$ and $I_v\subset G_v$ be its inertia subgroup.

  The following conditions on $\LL/\KK$ are introduced:

  \begin{description}
    \item[(LPTR)] \emph{Locally potentially totally ramified}: $I_v$ is open in $G_v$.
    \item[(NC)\ \ ] \emph{Normal closure}: For a Lie filtration $(\GG_n)$ of $\GG$, and the associated local filtration $(G_n)$ of $G$, the normal closure in $\GG/\GG_n$ of $G_{n-1}/G_n$ coincides with $\GG_{n-1}/\GG_n$.
    \item[(CE)\ \ ] \emph{Central element}: $\exists\,\tau\in Z(\GG)$ acting on $\mug_{p^\infty}(\LL)$ via $\zeta\mapsto \zeta^g$ with $g\in\ZZ, g>1$.
  \end{description}
  ( $\mug_{p^\infty}(\LL)$ denotes the group $\mug_{p^\infty}\cap \LL$, where $\mug_{p^\infty}$ is the group of $p$-power roots of unity.)
  \begin{theorem}\cite[Theorem 1]{ContiTerracini2025}\label{teo:ContiTerracini}
    If $\LL/\KK$ satisfies \textup{(LPTR)}, \textup{(NC)}, and \textup{(CE)} at some $v$ dividing $p$, then $\LL$ has the Bogomolov property \B.
  \end{theorem}
Condition {(LPTR) } makes Sen's theorem applicable to the local
$p$-adic Lie extension after a finite extension of the local
base field, yielding the metric inequalities described
in \S~4.1. Condition {(NC) } ensures that an element of
$\LL^{\GG_n}$ fixed by all conjugates of the local subgroup
$G_{n-1}/G_n$ already belongs to $\LL^{\GG_{n-1}}$.
Thus an element that does not belong to the preceding level
must be moved by one of these conjugate subgroups, making
the local metric inequality nontrivial at a suitable place
above $v$. Finally, condition {(CE)}, although somewhat ad hoc,
enables a descent argument used to deal with elements
$\alpha\in\LL^{\GG_n}$ such that
$\alpha^p\in\LL^{\GG_{n-1}}$. It automatically holds when $\mug_{p^\infty}(\LL)$ is finite; when this is not the case, it is equivalent to requiring that
the image of $Z(\GG)$ under the $p$-adic cyclotomic
character be open in $\ZZ_p^\times$;
see \cite[Proposition 2.8]{ContiTerracini2025}.

Recently, Theorem~\ref{teo:ContiTerracini} has been generalized
to the case where the base field \(\KK\) is not necessarily
a number field, provided that it has finite ramification
index at \(v\) and every finite extension of \(\KK\) has
property \B; see
\cite[Theorem~3.2]{ContiPirasTerracini2026}.

\subsection{Working on the prime-to-$p$ part: supercuspidal modular forms}
\label{sect:supercuspidal} 
In the proof of Theorem \ref{teo:AmorosoTerracini}, the strong supersingular condition $a_p=0$ was used in dealing with the prime-to-$p$-part in order to deduce that a power of $\Frob_p$ is central, and then exploit \ADZ lemma. At present, no new ideas have emerged that would allow the proof to be extended to the full generality predicted by Conjecture \ref{conj:modularB}.

 Nevertheless, the ADZ lemma can be successfully applied to another important class of modular forms in order to handle the prime-to-\(p\) part of the Bogomolov property. Namely, one may consider modular forms that are \emph{supercuspidal} at \(p\), i.e. those whose local Weil--Deligne representation at \(p\) is irreducible.
By \cite[\S 2.2.1]{Tate1979}, irreducibility implies that the projective image of such a representation is finite, leading to the following

\begin{proposition}\cite[Proposition 5.2]{ContiPirasTerracini2026}
Let $f\in S_k(\Gamma_1(N))$ be an eigenform, and let
\[\widehat\rho_f^{(p)}:G_\QQ \longrightarrow \GL_2\Big(\prod_{\lambda \nmid p} \KK_{f,\lambda}\Big)\] 
be the prime-to-$p$ part of the adelic Galois representation associated to $f$.
Let $\pi=\bigotimes_v \pi_v$ be the automorphic representation of $\GL_2(\mathbb{A}_\QQ)$ associated to $f$ and assume that its $p$-component $\pi_{p}$ is supercuspidal.
Then the projective image of $G_p$ through $\widehat\rho_f^{(p)}$ is finite, hence $\widehat\rho_f^{(p)}$ has \B by Lemma \ref{lem:ADZrepr}.
\end{proposition}

On the other hand the study of the $p$-part\[\widehat\rho_{f,p}:G_\QQ \longrightarrow \GL_2( \KK_f\otimes\QQ_p),\] 
for $f$ supercuspidal at $p$, is generally involved; indeed its restriction to $G_p$  is a product of potentially crystalline representations which can differ from each other,  although they share the same inertial type; it does not seem easy to show conditions \textup{}{(LPTR), (NC)} and \textup{(CE)} in this case. A partial result proving property \B for  the whole adelic representation associated to some  $p$-supercuspidal forms is presented in \cite[Theorem A]{ContiPirasTerracini2026}.

\section{Representations with big local image}\label{sect:applications}
In the remainder of this article, we present some applications  of Theorem~\ref{teo:ContiTerracini} to \(p\)-adic Galois representations, without assuming that they arise from modular forms or from geometry.

  Indeed, the theorem can be applied to an extension $\LL=\KK(\rho)$, where  $\KK$ is a number field and $\rho$ is a
continuous $p$-adic representation of $G_\KK$.

A particularly simple situation arises when the image of the inertia at a place over $p$ is \emph{full}:

\begin{theorem}\cite[Theorem 3.3]{ContiTerracini2025}
Let $\rho:G_\KK\rightarrow \GL_d(\ZZ_p)$ be a continuous representation
such that the image of an inertia subgroup at a place $v$ over $p$
is open in $\GL_d(\ZZ_p)$. Then $\rho$ has property \B.
\end{theorem}

Indeed, an open inertia image is an open $p$-adic Lie subgroup
of $\GL_d(\ZZ_p)$, hence has the same dimension as
$\GL_d(\ZZ_p)$. Conditions \textup{(LPTR)} and \textup{(NC)}
follow readily. As for condition \textup{(CE)}, we sketch the proof.
Let $U$ be the image of the inertia subgroup $I_v$ in
$\GL_d(\ZZ_p)$. Since $U$ is open, the local correspondence
between $p$-adic analytic groups and their Lie algebras
(see \cite[Section 9.5]{Dixon1999}) implies that the commutator
subgroup $[U,U]$ has finite index in
$U\cap\mathrm{SL}_d(\ZZ_p)$.
Consequently, every continuous character
$\chi:U\to\ZZ_p^\times$ satisfies a relation
\(
\langle\chi\rangle^x=\langle\det\rangle^y
\)
for suitable $x,y\in\ZZ_p$, with $x\neq0$, where
$\langle\cdot\rangle$ denotes the projection onto the
principal units. 
Let $\LL=\KK(\rho)$. If $\mug_{p^\infty}(\LL)$ is finite, condition
\textup{(CE)} holds automatically. Otherwise, the
cyclotomic character factors through $\Gal(\LL/\KK)$ and has
infinite image on $U$, so the corresponding exponent
$y$ is nonzero. Since $U$ contains an open subgroup
of scalar matrices, the cyclotomic character has infinite image
on these matrices, yielding condition \textup{(CE)}.

\subsection{The two-dimensional case}
 Recall that a representation $\rho$ of a profinite group $\Omega$ is said to be \emph{strongly absolutely irreducible} if its restriction to every open subgroup of \(\Omega\) is absolutely irreducible. \\

Let $\rho:\Omega\rightarrow \GL_d(\ZZ_p)$ be a continuous representation
of a profinite group, and let $G=\rho(\Omega)$.
Of course, if $G$ contains an open subgroup of $\SL_d(\ZZ_p)$,
then $\rho$ is strongly absolutely irreducible.
When $d=2$, the converse also holds, by an analysis of
the possible Lie algebras attached to $G$
(see \cite[Corollary B.1.3(i)]{ContiTerracini2025}).
For $d>2$, however, the converse fails: for instance,
$\mathrm{Sym}^{d-1}\SL_2(\ZZ_p)$ acts strongly absolutely
irreducibly on $\QQ_p^d$, but does not contain any open
subgroup of $\SL_d(\ZZ_p)$, since its dimension is
$3<d^2-1$.

Applying the above considerations to a 2-dimensional $p$-adic Galois representation leads to the following result:

\begin{proposition}
 Let $\rho : G_\KK \rightarrow \GL_2(\ZZ_p)$ be a continuous representation. 
 Let $v$ be a $p$-adic place of $\KK$, with
decomposition group $G_v$ and inertia subgroup $I_v$.\\
Assume that $\rho|_{G_v}$ is strongly absolutely irreducible, and
that $\det \rho(I_v)$ is infinite. \\Then $\rho(I_v) $ is open in $\GL_2(\ZZ_p)$, and hence $\rho$ has property \B.
\end{proposition}

Indeed, strong absolute irreducibility implies that
$\rho(G_v)$ contains an open subgroup of $\SL_2(\ZZ_p)$.
Since $G_v/I_v$ is abelian, $\rho(I_v)$ contains the
commutator subgroup of $\rho(G_v)$, and hence also an open
subgroup of $\SL_2(\ZZ_p)$.
The hypothesis that $\det\rho(I_v)$ is infinite implies
that it is open in $\ZZ_p^\times$. Together, these two
properties imply that $\rho(I_v)$ is open in
$\GL_2(\ZZ_p)$.
For details, see \cite[Proposition B.2.10]{ContiTerracini2025}.

\paragraph{A criterion for strong absolute irreducibility in dimension \(2\).}
When \(d=2\), strong absolute irreducibility can be detected by means of the following characterization \cite[Corollary B.1.3]{ContiTerracini2025}:

\begin{proposition}\label{prop:criterioSI}
A two-dimensional representation  of a profinite group with values in $\QQ_p$ is strongly absolutely irreducible if and only if
it  is absolutely irreducible, not induced from a character, and not a twist of a finite-image representation.
\end{proposition}

\subsection{Applications}
In this section, we present some applications of the criteria discussed above which, in our opinion, are particularly noteworthy. The interested reader is referred to \cite[Sections 4--7]{ContiTerracini2025} for a more extensive exposition.
\paragraph{The $p$-part of supersingular modular forms.}

We recall that a normalized cuspidal eigenform
\[
f=\sum_{n\geq 1} a_n q^n \in S_k(\Gamma_1(N))
\]
is said to be \emph{supersingular} at \(p\) if $p$ does not divide $N$ and \(v_p(a_p)>0\). In this case the local $p$-adic representation is  crystalline.   When $k=2$ and $a_p\in\ZZ$, as is the case for elliptic
curves over $\QQ$, supersingularity and the Hasse bound
force $a_p=0$ for $p\geq5$; this is the situation referred
to as \emph{strong supersingularity}. For general weights and Hecke fields \(\KK_f\), supersingularity no longer implies strong supersingularity, and forms with \(v_p(a_p)>0\) can satisfy \(a_p\neq 0\).

 The classification of irreducible two-dimensional crystalline representations with coefficients in $\QQ_p$ given in \cite[Prop. 3.1.1]{Breuil2003}, combined with Proposition \ref{prop:criterioSI}, leads to the following

\begin{theorem}\label{teo:apnotzero} Let $f=\sum_{n\geq 1} a_n q^n \in S_k(\Gamma_1(N))$ be a normalized Hecke eigenform of weight $k\geq 2$, with Hecke field $\KK_f$. Let $p$ be a prime not dividing $N$ and $v$ a splitting $p$-adic place of $\KK_f$ such that $f$ is supersingular at $v$.
 If $a_p\not=0$ then $\rho_{f,v}$ has \B.
\end{theorem}

It is worth noting that this criterion points in a direction opposite to that of Habegger's and Amoroso--Terracini's arguments: the strong supersingularity condition \(a_p=0\) occurs precisely when the local representation is induced from a character, a situation excluded by strong absolute irreducibility.\\

For example, the \emph{Ramanujan $\Delta$ function}
(labelled \texttt{1.12.a.a} in \cite{LMFDB}) is an eigenform
in $S_{12}(\Gamma_0(1))$ with rational Hecke field and
$a_7=-16744\neq0$, with $a_7\equiv0\pmod7$.
It follows from Theorem~\ref{teo:apnotzero} that
$\rho_{\Delta,7}$ has property \B.

An example with $[\KK_f:\QQ]>1$ is the eigenform
$f\in S_2(\Gamma_0(43))$ labelled \texttt{43.2.a.b}
in \cite{LMFDB}. Here $\KK_f=\QQ(\sqrt2)$ and the  minimal polynomial of $a_{17}$ over $\QQ$
is $X^2-10X+17$.
The prime $17$ splits in $\KK_f$, and $a_{17}$ has valuation
$1$ at one of the two places above $17$ and valuation $0$
at the other. Choosing the former place $v$, we obtain
$\KK_{f,v}\simeq\QQ_{17}$ and $v(a_{17})>0$.
Thus Theorem~\ref{teo:apnotzero} applies, showing that
$\rho_{f,v}$ has property \B. 
 
\paragraph{Elliptic curves over a number field.}

Habegger's argument does not extend to elliptic curves defined over
an arbitrary number field, essentially because Elkies' theorem is not known in
such generality. In fact, it is proved in \cite[Theorem 1.4]{Sahu2025} that
the Bogomolov property holds for the adelic representation
\[
\widehat{\rho}_E:G_\KK \longrightarrow \GL_2(\widehat{\ZZ})
\]
associated with an elliptic curve \(E\) defined over a number field \(\KK\),
provided that the curve has infinitely many supersingular primes. In particular,
this applies to every elliptic curve when \(\KK\) has at least one real
embedding \cite[Corollary 1.5]{Sahu2025}.

Outside this setting, Theorem \ref{teo:ContiTerracini} sometimes allows us to
recover property \B for the \(p\)-adic representation
\[
\rho_{E,p}:G_\KK\longrightarrow \GL_2(\ZZ_p).
\] or, equivalently, for the field \(\KK(E[p^\infty])\). 

\begin{theorem}\cite[Theorem 6.1]{ContiTerracini2025}
Assume that \(E_{/\KK}\) is non-CM and that there exists a rational prime \(p\)
such that:
\begin{itemize}
\item[i)] for some \(p\)-adic place \(v\) of \(\KK\), the restriction
\(
\rho_{E,p}|_{G_v}
\)
is absolutely irreducible;
\item[ii)] the residual representation
\(
\bar{\rho}_{E,p}:G_\KK\longrightarrow \GL_2(\FF_p)
\)
is surjective.
\end{itemize}
Then \(\rho_{E,p}\) has property
\B.
\end{theorem}

By Serre's big image theorem, for any fixed non-CM elliptic curve \(E\), condition~$ii)$ holds for all sufficiently large primes \(p\).

\paragraph{Property \B in families.}
One of the advantages of the methods described above is that they behave well in \(p\)-adic analytic families of Galois representations.

This framework encompasses several well-known examples, including:
\begin{itemize}
\item universal deformation spaces attached to a fixed residual representation \(\bar\rho\);
\item Hida families;
\item non-ordinary components of the Coleman--Mazur eigencurve.
\end{itemize}

Let
\[
\{\rho_x:G_\KK\to\GL_2(\ZZ_p)\}_{x\in X}
\]
be such a family over a rigid-analytic space \(X\), and fix a place \(v\) of $\KK$, dividing $p$.

After removing a finite union of proper analytic subsets, and under suitable hypotheses, one obtains an open subset \(U\subset X\) such that for every \(x\in U\):
\begin{enumerate}
\item \(\rho_x|_{G_v}\) is absolutely irreducible and not induced;
\item the Hodge--Tate--Sen weights satisfy \(k_{x,1}\neq k_{x,2}\);
\item \(\det\rho_x(I_v)\) has infinite image in \(\ZZ_p^\times\).
\end{enumerate}
By applying Theorem \ref{teo:ContiTerracini}, the following 
\textit{generic Bogomolov property} is deduced:
\begin{center}\textit{For every \(x\in U\), the representation \(\rho_x\) has property \B.}\end{center}

The interested reader can find precise statements and complete proofs in Section~7 of \cite{ContiTerracini2025}.

A natural question concerns the existence of a uniform lower bound for the height throughout the family. We hope to return to this problem in future work.

\paragraph{Beyond condition (CE): the $p$-adic Rémond conjecture.}
Although somewhat ad hoc and rather mysterious, the central
element condition plays a crucial role in the proof of
Theorem~\ref{teo:ContiTerracini}.
Consider, for example, the extension $\LL/\QQ$, with
\[
\LL=\QQ\bigl(\mug_{p^\infty},\,\sqrt[p^n]{2}\,:\,n\in\NN\bigr).
\]
It is a $p$-adic Lie extension whose Galois group $\GG$
is an open subgroup of $\ZZ_p^\times\ltimes\ZZ_p$,
equal to the whole group when $p$ is odd.
Denote by
$G_p\subseteq\GG$ and $I_p\subseteq G_p$ the decomposition
and inertia subgroups of $\LL$ at a place $v$ above $p$.

Condition (LPTR) is satisfied: indeed, $\LL/\QQ$ is totally
ramified at $2$ when $p=2$; for odd $p$, it is totally
ramified at $p$ over a finite intermediate extension
(see \cite[Theorems 5.4 and 5.5]{Viviani2004}).
For odd $p$, the same references also imply that $I_p$
is open in $\GG$; for $p=2$, we have $I_p=\GG$.
Let $(\GG_n)$ be a Lie filtration of $\GG$, and put
$G_n=G_p\cap\GG_n$.
For all sufficiently large $n$, we have $\GG_n\subseteq I_p$,
and hence $G_n=\GG_n$.
Thus $G_{n-1}/G_n=\GG_{n-1}/\GG_n$, so condition (NC)
holds after shifting the indices of the filtration.
However, condition (CE) fails, because the center of $\GG$
is trivial and $\LL$ contains all $p$-power roots of unity.
The field $\LL$ does not have property \B, since
\[
h\bigl(\sqrt[p^n]{2}\bigr)=\frac{\log 2}{p^n}
\longrightarrow 0 \qquad\text{as }n\longrightarrow\infty.
\]

Other sequences in $\LL$ whose heights tend to zero can
of course be obtained from $(\sqrt[p^n]{2})$ by multiplying
its terms by roots of unity. A natural question is whether
all sequences in $\LL^\times$ whose heights tend to zero
are essentially of this form, i.e., must eventually lie in the
divisible hull of the subgroup generated by $2$.\\
Recall that the central element condition was introduced
to overcome the obstruction to descent posed by successive
$p$-th roots, which can give rise to sequences of elements
whose heights tend to zero. Kummer theory detects radical
elements in finite Galois extensions through the condition
that their conjugates differ from them by multiplication
by roots of unity. This suggests that, even when condition (CE) fails, our  height argument
might still yield a positive lower bound outside a suitable
exceptional set of radical elements.

The above considerations fit in the broader context of the following conjecture, due to Rémond \cite[Conjecture 3.4]{Remond2017}:

\begin{conjecture}[Rémond]\label{conj:Remond}
Let $\KK$ be a number field, and 
$\Gamma \subseteq  \KK^\times$ be a non-torsion finitely generated subgroup.
Put $\LL=\KK(\Gamma^{\rm{div}})$.
 The Weil height is bounded from below by a strictly positive constant on $\LL^\times \setminus \Gamma^{\rm{div}}$.
\end{conjecture}

Here $\Gamma^{\rm{div}}$ denotes the divisible subgroup of $\Qbar^\times$ generated by $\Gamma$:
\[\Gamma^{\rm{div}}=\{\alpha\in\Qbar^\times\ |\ \alpha^n\in \Gamma, \hbox{ for some } n\in\NN\}.\]
Although Rémond's conjecture appears to lie beyond the
scope of our methods, its $p$-adic counterpart can be
approached using the techniques on $p$-adic Lie extensions
described above. This problem has been addressed in
\cite{ContiDelCorsoPlessisTerracini2026}.

It turns out that Conjecture~\ref{conj:Remond} holds if we
restrict our attention to the field generated by the
$p$-divisible hull of $\Gamma$, giving the following
$p$-adic version:

 \begin{theorem}[\cite{ContiDelCorsoPlessisTerracini2026}]\label{teo:Remond-p-adic}
 Let $\KK$ be a number field, and 
$\Gamma \subseteq \KK^\times$ be a non-torsion finitely generated subgroup.
Put $\LL=\KK(\Gamma^{p-\rm{div}})$.
 The Weil height is bounded from below by a strictly positive constant on $\LL^\times \setminus \Gamma^{\rm{div}}$.
\end{theorem}

Here $\Gamma^{p-\rm{div}}$ denotes the $p$-divisible subgroup of $\Qbar^\times$ generated by $\Gamma$:
\[\Gamma^{p-\rm{div}}=\{\alpha\in\Qbar^\times\ |\ \alpha^{p^n}\in \Gamma, \hbox{ for some } n\in\NN\}.\]

The exceptional set to be excluded must involve the full divisible hull
$\Gamma^{\rm div}$; indeed, the field
$\KK(\Gamma^{p\text{-}\rm div})$ may contain sequences
of elements of arbitrarily small positive height outside
$\Gamma^{p\text{-}\rm div}$, arising for example from roots
of unity or other prime-to-$p$ radicals already
contained in the base field.

When $\mathrm{rk}(\Gamma)=1$,  Theorem \ref{teo:Remond-p-adic} was proved by  Amoroso \cite{Amoroso2016}, for $\KK=\QQ,\, \Gamma=\langle \alpha\rangle$ where $\alpha\in \ZZ$ and $p$ is an odd rational prime not satisfying the Wieferich condition (that is, $p$ does not divide $\alpha$ and $p^2$ does not divide $\alpha^{p-1}-1$). 
A generalization is due to Plessis \cite[Théorème 1.8]{Plessis2022}, who established the same result when $\KK$ is any number field, $\alpha$ is any element of $\KK$, and $p$ is any odd rational prime not dividing the discriminant of $\KK$ . 

For the proof of the general case of Theorem~\ref{teo:Remond-p-adic}, the authors of \cite{ContiDelCorsoPlessisTerracini2026} used a modified version of Theorem~\ref{teo:ContiTerracini}, which does not require condition \textup{(CE)}, yielding a positive lower bound for the height outside an exceptional set; unlike the rank-one case, the verification of conditions (LPTR) and  \textup{(NC)} is no longer immediate, but can nevertheless be carried out in full generality. 

Moreover, they addressed the problem of dealing with a finite set $S$ of primes instead of a single prime $p$. This approach uses a recursive argument on $|S|$, based on the possibility of bounding below the height on composita of fields, under some conditions \cite[Proposition 4.2] {ContiPirasTerracini2026}. The  result they obtained, which contains Theorem \ref{teo:Remond-p-adic} as a particular case, is the following: 

\begin{theorem} (\cite[Corollary 0.7]{ContiDelCorsoPlessisTerracini2026})
Let $\Gamma\subseteq \overline{\QQ}^\times$ be a finitely generated subgroup, and let $S$ be a finite set of rational primes.
Assume that for all but at most one prime $p\in S$, there is a field embedding $\iota_p : \Qbar\hookrightarrow \overline{\QQ_p}$ such that $\iota_p(\Gamma)\subseteq \overline{\ZZ_p}^\times$. 
Put
\[
\Gamma^{S-\mathrm{div}}
=
\{\alpha\in\Qbar^\times :
\alpha^m\in\Gamma \text{ for some integer } m\geq 1
\text{ whose prime divisors belong to } S\}.
\]
Then the Weil height is bounded from below by a strictly positive constant on  the set $\KK(\mug_\infty, \Gamma^{S-\mathrm{div}})^\times\setminus \Gamma^{\mathrm{div}}$, for every finite extension $\KK/\QQ(\Gamma)$. 
\end{theorem}

\vspace{0.5cm} \indent {\it
A\,c\,k\,n\,o\,w\,l\,e\,d\,g\,m\,e\,n\,t\,s.\;} I would like to warmly thank the referee for the care taken in reading the manuscript, the interest in my work, and the helpful suggestions. I feel that these have enriched the paper and improved the new version. Their comments gave me the opportunity to reflect further on several aspects of the paper and broaden my perspective on the subject.\\
I am also grateful to Andrea Conti for his support and valuable comments.

\bigskip




\bigskip
\bigskip
\begin{minipage}[t]{10cm}
\begin{flushleft}
\small{
\textsc{Lea Terracini}
\\Dipartimento di Informatica
\\Università di Torino,
\\Corso Svizzera 185
\\10142, Torino, Italy
\\*e-mail: lea.terracini@unito.it

}
\end{flushleft}
\end{minipage}

\end{document}